\documentclass[conferences]{IEEEtran}
\usepackage{color}
\usepackage[pdftex]{graphicx}
\graphicspath{{pics/}}
\usepackage[cmex10]{amsmath}
\usepackage{amssymb,amsthm}
\usepackage{algpseudocode}
\usepackage[linesnumbered,ruled,vlined]{algorithm2e}
\usepackage{stmaryrd}
\usepackage{fixltx2e}
\usepackage{url}
\usepackage{multirow,booktabs}
\usepackage[caption=false,font=footnotesize,labelformat=simple]{subfig}

\newcommand\figref{Fig.~\ref}
\DeclareMathOperator*{\argmin}{argmin}
\newcommand\norm[1]{\left\lVert#1\right\rVert}
\newcommand\abs[1]{\lvert#1\rvert}

\usepackage{tabularx}
\makeatletter
\def\hlinewd#1{%
\noalign{\ifnum0=`}\fi\hrule \@height #1 %
\futurelet\reserved@a\@xhline}
\makeatother

\usepackage{threeparttable}

\usepackage{cite}
\ifCLASSINFOpdf
\else
\fi

\begin{document}
\title{Hierarchical Distributed EV Charging Scheduling in Distribution Grids}

\author{\IEEEauthorblockN{Behnam Khaki\IEEEauthorrefmark{1}, Yu-Wei Chung\IEEEauthorrefmark{2}, Chicheng Chu\IEEEauthorrefmark{3}, and Rajit Gadh\IEEEauthorrefmark{4}}\\
\IEEEauthorblockA{Smart Grid Energy Research Center (SMERC),
University of California, Los Angeles\\
Los Angeles, CA, USA, 90095}
{Email:\{ \IEEEauthorrefmark{1}behnamkhaki,
\IEEEauthorrefmark{2}ywchung,
\IEEEauthorrefmark{3}peterchu,
\IEEEauthorrefmark{4}gadh}\}@ucla.edu}


%Telephone: (800) 555--1212, Fax: (888) 555--1212}


\maketitle

\begin{abstract}
In this paper, a hierarchical distributed method consisting of two iterative procedures is proposed for optimal electric vehicle charging scheduling (EVCS) in the distribution grids. In the proposed method, the distribution system operator (DSO) aims at reducing the grid loss while satisfying the power flow constraints. This is achieved by a consensus-based iterative procedure with the EV aggregators (Aggs) located in the grid buses. The goal of aggregators, which are equipped with the battery energy storage (BES), is to reduce their electricity cost by optimal control of BES and EVs. As Aggs' optimization problem increases dimensionally by increasing the number of EVs, they solved their problem through another iterative procedure with their customers. This procedure is implementable by exploiting the mathematical properties of the problem and rewriting Aggs' optimization problem as the \textit{sharing problem}, which is solved efficiently by the alternating direction method of multipliers (ADMM). To validate the performance, the proposed method is applied to IEEE-13 bus system. 
\end{abstract}

% \fi
\IEEEpeerreviewmaketitle

\section{Introduction}
Transportation sector consumes a significant percentage of energy, and it has a considerable contribution to the air pollution and greenhouse gas emission. Over the last decade, it has been shown that electric vehicles (EVs) are a promising technology to reduce transportation's dependency on fossil fuels. However, due to EVs' electrical energy demand, they introduce new challenges to the electricity sector. EV charging load demand, in high penetration scenarios which is feasible in the near future, may lead to stability, quality, and economy issues in power grids. According to the EV load characteristics, it is considered as a controllable load which its adverse effects can be mitigated through a demand management strategy. Nonetheless, due to the uncertainty in EV load demand \cite{Behnam1}-\cite{Behnam2} and scalability issue in the case of high EV penetration \cite{Behnam3}, EV load management is challenging. In this paper, the scalability issue is addressed by introducing a hierarchical and fully distributed EV charging scheduling (EVCS). 

There is a rich body of literature proposing either centralized or distributed EVCS methods. In centralized approaches \cite{Clement-Nyns1}-\cite{Tang1}, the distribution system operator (DSO) receives (or predicts) the data relating to arrival time, departure time, and required energy of each individual EV, and it coordinates their charging load considering the defined objective function. Centralized methods, however, suffer from curse of dimensionality issue if DSO has to deal with a large population of EVs. In addition, centralized methods can not preserve the EV owners' privacy, as they have to communicate their sensitive information with DSO. To tackle these issues, researchers propose distributed methods in which DSO solves EVCS problem through an iterative procedure with EV aggregators (Aggs) \cite{Mohiti}, EVs \cite{Callaway}, or both \cite{Xu}-\cite{Giannakis}. In the first case, it is Agg who directly receives EVs' information and is in contact with DSO. In the second case, DSO executes the iterative procedure directly with EVs, so they do not need to share their sensitive information with any entity. In the last case, DSO communicates with Aggs, and each Agg communicates with its EVs, and no sensitive information is shared with other entities. Obviously, the last case is more scalable, and its structure has more flexibility from different entities' objective function perspective. It is worthwhile to mention that among the distributed methods, some of them do not consider the power grids model, while others take the power flow constraints into consideration \cite{Behnam3}-\cite{Amini}. In this paper, the focus is only on the second group.

Among the proposed distributed methods, the authors in \cite{Mohiti} use alternating direction method of multipliers (ADMM) to solve EVCS problem. The authors in \cite{Callaway} use a novel shrunken primal-dual subgradient method for valley-filling problem. The authors in \cite{Xu} benefit from \cite{Low5} to solve Agg's problem by a hierarchical and distributed method. In \cite{Giannakis}, the authors use Frank-Wolfe method to make the EVCS problem scalable, where the optimization problem is formulated as a linear program. Nevertheless, the proposed methods are developed either based on strict assumptions or for specific objective functions. To address those issues, we further extend our previous work \cite{Behnam3} and propose a fully distributed EVCS by the consideration of power flow constraints. In our method which is based on multi-agents systems, DSO, Aggs, and EVs solve their objective function locally through a hierarchical iterative communication procedure implemented by ADMM.

The paper is structured as follows: in Section \ref{System Model Description}, the power flow constraints as well as EV and battery energy storage (BES) models are introduced; in Section \ref{Problem Formulation}, the EVCS problem is formulated, and it is solved by our hierarchical distributed method based on ADMM; numerical simulation results of the proposed EVCS are shown in Section \ref{Simulation Results}; and the paper is concluded in Section \ref{Conclusion}.  
      
\section{System Model Description}  \label{System Model Description}
\subsection{Distribution Grid Model}
\begin{figure}[h]%\hspace*{-4cm}
    \centering \includegraphics[clip,scale=0.9, trim=9.7cm 10.9cm 12.5cm 8.6cm] %[clip,scale=0.38, trim=9.0cm 3.0cm 3.5cm 1.85cm] 
    {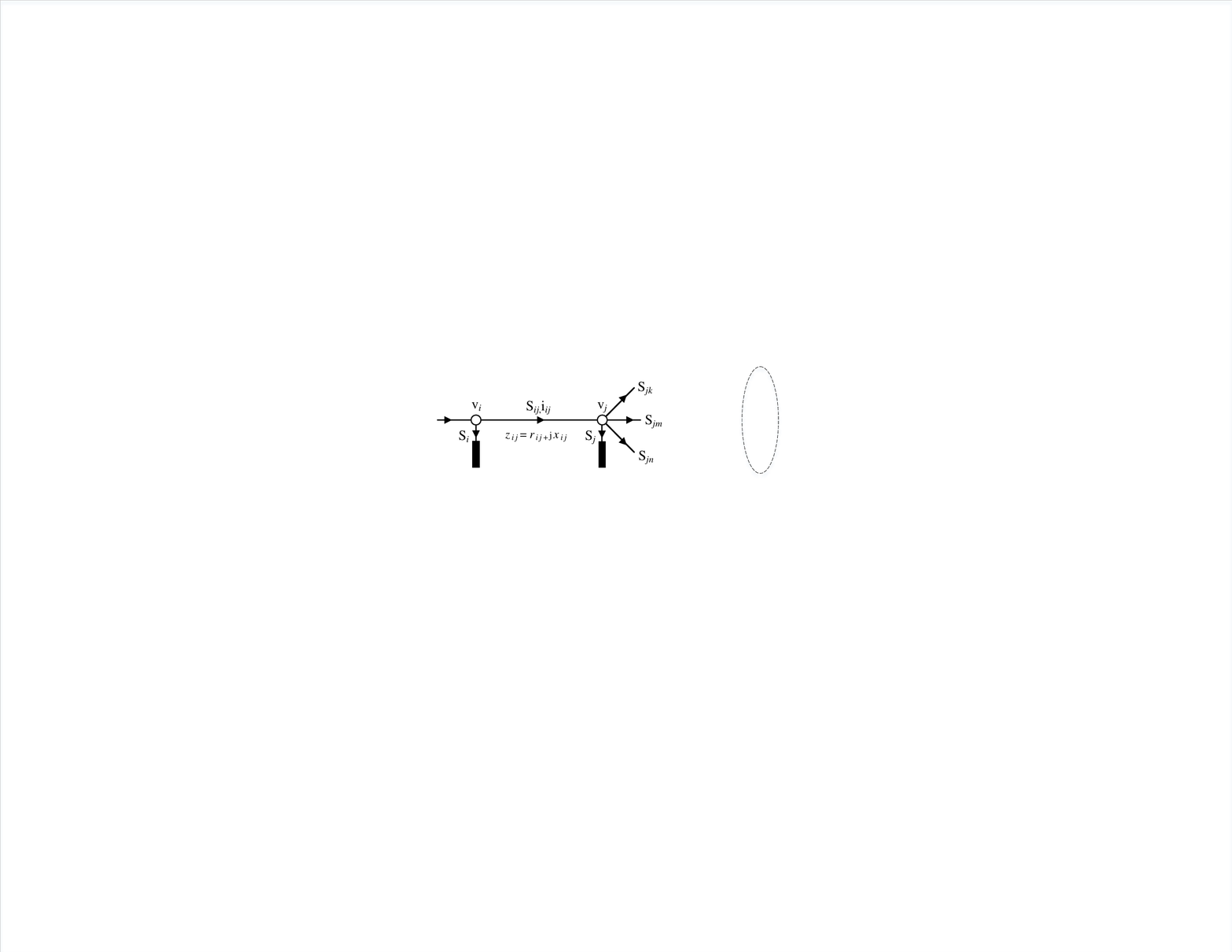}%{Grid3.pdf}%{Grid2.jpg}%{Fig1.jpg} %{PhysicalLayer.jpg}%
    \caption{DistFlow model of the distribution line.}
    \label{fig:DistFlow}
\end{figure}

In this paper, we consider the distribution grid modeled as a connected graph which is shown by $G=(\mathbb{N}_b,\zeta)$, where $\mathbb{N}_b$ denotes the set of the grid buses, and $\zeta$ denotes the set of the lines. we use $\mathbf{y}_{n}:=\big(y_{n}(t),y_{n}(t+1), \dots,y_{n}(t+N-1)\big)^{T}$, where $N \in \mathbb{N}$ is the time horizon, for all the variables. If $(i,j) \in \zeta$, a line connects bus $i$ to bus $j$, which its impedance and current are shown by $z_{ij}=r_{ij}+\text{j}x_{ij}$, $z_{ij} \in \mathbb{Z}$, and $\mathbf{i}_{ij}$, respectively. The apparent power of the line $(i,j) \in \zeta $ is shown by $\mathbf{S}_{ij}$, $\mathbf{S}_{ij} \in \mathbb{Z}^{N}$, and it is calculated by $\mathbf{S}_{ij}=\mathbf{P}_{ij}+\text{j}\mathbf{Q}_{ij}$, where $\mathbf{P}_{ij}$ and $\mathbf{Q}_{ij}$ are the active and reactive powers, respectively. Also, the voltage, active power, reactive power, and apparent power at bus $i$ are shown by $\mathbf{v}_i$, $\mathbf{P}_{i}$, $\mathbf{Q}_{i}$, and $\mathbf{S}_{i}$, respectively. To show the relations between bus voltages, active powers, reactive powers, line currents and impedances, DistFlow model \cite{Baran1} is assumed (\figref{fig:DistFlow}) for the distribution grid as follows:
\begin{subequations} \label{Grid1}
\begin{align}
&\mathbf{v}_i-\mathbf{v}_{j}=z_{ij}\mathbf{i}_{ij} \label{Grid11}\\
&\mathbf{S}_{ij}=\mathbf{v}_i {\mathbf{i}_{ij}}^* \label{Grid12}\\
&\mathbf{S}_{j}=\mathbf{S}_{ij}-z_{ij}{\abs{\mathbf{i}_{ij}}}^2-\sum\limits_{(j,k) \in \zeta}^{}{\mathbf{S}_{jk}}. \label{Grid13}
\end{align}
\end{subequations}

Substituting \eqref{Grid11} and \eqref{Grid12} in \eqref{Grid13}, the following equations are obtained,
\begin{subequations} \label{Grid2}
\begin{align}
&\mathbf{P}_{j}=\mathbf{P}_{ij}-r_{ij}\mathbf{I}_{ij}-\sum\limits_{(j,k) \in \zeta}^{}{\mathbf{P}_{jk}} \label{Grid21}\\
&\mathbf{Q}_{j}=\mathbf{Q}_{ij}-x_{ij}\mathbf{I}-\sum\limits_{(j,k) \in \zeta}^{}{\mathbf{Q}_{jk}} \label{Grid22}\\
&{\mathbf{V}}_i-{\mathbf{V}}_j = 2(r_{ij}\mathbf{P}_{ij}+x_{ij}\mathbf{Q}_{ij})-({r_{ij}}^2+{x_{ij}}^2)\mathbf{I}_{ij} \label{Grid23}\\
&{\mathbf{V}}_i{\mathbf{I}}_{ij}=\mathbf{P}_{ij}^2+\mathbf{Q}_{ij}^2\label{Grid24},
\end{align}
\end{subequations}
where $\mathbf{V}_i={\abs{\mathbf{v}}_i}^2$ and $\mathbf{I}_i={\abs{\mathbf{i}}_i}^2$. The voltage at the root node ($\mathbf{v}_1$) should be equal to a constant value, assuming the distribution grid is connected to an infinite bus, while the other voltages can vary within a limited range. This is shown by:
 \begin{subequations} \label{Grid3}
\begin{align}
&\underline{\mathbf{v}}\leq{\mathbf{v}_i}\leq{{\overline{\mathbf{v}}}},~i \in \mathbb{N}_b\setminus{\{1\}}\label{Grid31}\\ 
& {\mathbf{v}}_1 = \mathbf{v}_{ref}. \label{Grid32}
\end{align}
\end{subequations}

\subsection{EV Charging Model}  \label{EV Charging Model}
The set of EVs supplied by $Agg_j$ and their number are shown by $\mathbb{N}_{j,ev}$ and $\mathcal{N}_{j,ev}$, respectively. Also, $EV_{j,i}$ stands for \textit{i}\textsubscript{th} EV supplied by $Agg_j$. It is assumed that each $EV_{j,i}$ is located in either a commercial or a residential building. We use $EVB_{j,i}$ for each $EV_{j,i}$ and its corresponding building, which is modeled as a discrete-time linear system as follows \cite{Behnam3}: 
\begin{subequations} \label{RESModel1}
\begin{align}
&c_{j,i}^{ev}(t+1)=c_{j,i}^{ev}(t)+T_h\eta_{j,i}^{ev} p_{j,i}^{ev}(t) \label{RESModel11}\\
&e_{j,i}^{evb}(t)=p_{j,i}^{uc}(t)+p_{j,i}^{ev}(t), \label{RESModel12}
\end{align}
\end{subequations}
where $c_{j,i}^{ev}, p_{j,i}^{ev}, p_{j,i}^{uc}, e_{j,i}^{evb} \in \mathbb{R}$, $\eta_{j,i}^{ev} \in \mathbb{R}^{+}_{\leqslant{1}}$. $c_{j,i}^{ev}(t)$ is the energy stored in EV battery at time $t$, $p_{j,i}^{uc}(t)$ is the non-EV and uncontrollable active load demand minus the power generated by the solar panel of the building, $T_h$ is $0.5$ in this paper which corresponds to $30$ min, and $p_{j,i}^{ev}(t)$ is the control variable which is determined by the EVCS. As we assume that each $EVB_{j,i}$ is provided by a solar panel and has an EV charger with vehicle-to-grid (V2G) capability, it may supply power to the grid.    

The constraints on the EV charging/discharging, relating to the charger power rating and the EV battery capacity, are:
\begin{subequations} \label{RESModel2}
\begin{align}
&\underline{\mathbf{p}}_{j,i}^{ev}\leqslant{\mathbf{p}}_{j,i}^{ev}\leqslant{\overline{\mathbf{p}}_{j,i}^{ev}} \label{RESModel21}\\
&\underline{\mathbf{C}}_{j,i}^{ev}\leqslant{\mathbf{c}_{j,i}^{ev}}\leqslant{\overline{\mathbf{C}}_{j,i}^{ev}}, \label{RESModel22}
\end{align}
\end{subequations}
where $\underline{C}_{j,i}^{ev}(t), \overline{C}_{j,i}^{ev}(t) \in \mathbb{R}$ are the EV battery time-varying constraints which are defined as follows; if $EV_{j,i}$ is: 
\begin{itemize}
\item{not plugged in $EVB_{j,i}$, $\underline{C}_{j,i}^{ev}(t)=\overline{C}_{j,i}^{ev}(t)=0$.} 
\item{plugged in $EVB_{j,i}$, but it is in idle mode, ${\underline{C}_{j,i}^{ev}(t)=0}$ \& $\overline{C}_{j,i}^{ev}(t)=C_{j,i}^{ev}$, where $C_{j,i}^{ev} \in \mathbb{R}$ is the maximum EV battery energy capacity.} 
\item{plugged in $EVB_{j,i}$, and it is needed by time $t$, $\underline{C}_{j,i}^{ev}(t)=\overline{C}_{j,i}^{ev}(t)=C_{j,i}^{ev}$.} 
\end{itemize}

We define the set of feasible charging trajectories of $EV_{j,i}$ as:
%of $EVC_{j,i}$, $i \in \mathbb{N}_{\mathcal{I}}$, which is supplied by $FA_j$, $j \in \mathbb{N}_{\mathcal{J}}$, as: 
\begin{align}
\label{RESSet}
\mathbb{U}_{j,i}^{ev} =\bigg\{{\mathbf{p}_{j,i}^{ev} \in \mathbb{R}^{N}}\mid{\eqref{RESModel1}-\eqref{RESModel2}~\forall t \in \llbracket k,k+N-1\rrbracket} \bigg\}.
\end{align}

\subsection{BES Model}  \label{BES Model}
The BES controlled by $Agg_j$ is indicated by $BES_j$ and modeled as follows:
\begin{align}
\label{BES}
&c_{j}^{bes}(t+1)=c_{j}^{bes}(t)+T_h\eta_{j}^{bes} p_{j}^{bes}(t),
\end{align}
where $c_{j}^{bes}$ and $p_{j}^{bes} \in \mathbb{R}$, $\eta_{j}^{bes} \in \mathbb{R}^{+}_{\leqslant{1}}$. The energy stored in BES and its charging/discharging power are limited by the following constraints:
\begin{subequations} \label{BESConstraint}
\begin{align}
&{\mathbf{s}_{j}^{bes}}^2={\mathbf{p}_{j}^{bes}}^2+{\mathbf{q}_{j}^{bes}}^2\leqslant{{\overline{\mathbf{s}}_{j}^{bes}}^2} \label{BESConstraint2}\\
&\underline{\mathbf{C}}_{j}^{bes}\leqslant{\mathbf{c}_{j}^{bes}}\leqslant{\overline{\mathbf{C}}_{j}^{bes}}, \label{BESConstraint1}
\end{align}
\end{subequations}
where $\mathbf{p}_{j}^{bes}$, $\mathbf{q}_{j}^{bes}$ and $\mathbf{s}_{j}^{bes}$ are the active, reactive, and apparent powers, respectively, $\overline{\mathbf{s}}_{j}^{bes}$ is the apparent power rating of the $BES_j$'s bi-directional converter, and $\underline{C_{j}}^{bes}$ and $\overline{C_{j}}^{bes} \in \mathbb{R}$ are the minimum and maximum $BES_j$ energy constraints, respectively.
\subsection{EV Aggregator Model}  \label{Aggregator Model}
We use $\mathbb{N}_{Agg}$ to denote the set of buses which have EV aggregator. According to the model defined for EVBs, the active and reactive powers of bus $j$, where $Agg_j$ is located, are obtained as:
\begin{subequations} \label{Grid4}
\begin{align}
&\mathbf{P}_{j_c}=\mathbf{p}_{j}^{bes}+\sum\limits_{i \in \mathbb{N}_{j,ev}}^{}{(\mathbf{p}_{j,i}^{ev}+\mathbf{p}_{j,i}^{uc})}\label{Grid41}\\
&\mathbf{Q}_{j_c}=\mathbf{q}_{j}^{bes}+\sum\limits_{i \in \mathbb{N}_{j,ev}}^{}{\mathbf{q}_{j,i}^{uc}}, \label{Grid42}
\end{align} 
\end{subequations}
in which, $\mathbf{q}_{j,i}^{uc}$ is the uncontrollable reactive load demand by $EVB_{j,i}$. As $\mathbf{p}_{j,i}^{uc}$ and $\mathbf{q}_{j,i}^{uc}$ are not controllable, we consider their aggregated value at $Agg_j$ in the EVCS modeling, which are defined as:
\begin{subequations} \label{Grid5}
\begin{align}
&\mathbf{p}_{j}^{uc}=\sum\limits_{i \in \mathbb{N}_{j,ev}}^{}{\mathbf{p}_{j,i}^{uc}} \label{Grid51}\\
&\mathbf{q}_{j}^{uc}=\sum\limits_{i \in \mathbb{N}_{j,ev}}^{}{\mathbf{q}_{j,i}^{uc}}. \label{Grid52}
\end{align} 
\end{subequations} 

\begin{figure}[h]%\hspace*{-4cm}
    \centering \includegraphics[clip,scale=0.5, trim=6.0cm 7.5cm 7.5cm 4.0cm] %[clip,scale=0.38, trim=9.0cm 3.0cm 3.5cm 1.85cm] 
    {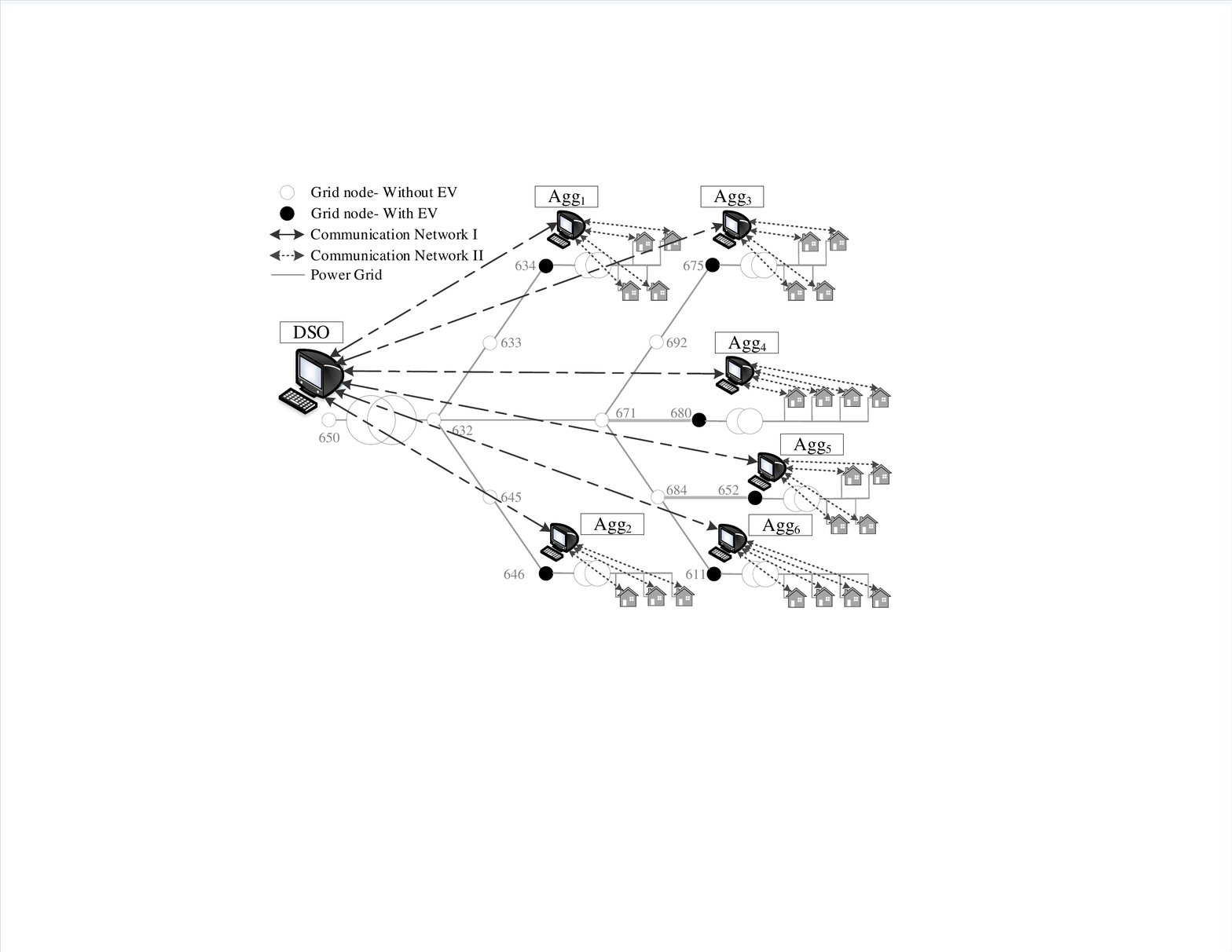}%{Grid3.pdf}%{Grid2.jpg}%{Fig1.jpg} %{PhysicalLayer.jpg}%
    \caption{IEEE-13 bus system with the proposed hierarchical distributed EVCS.}
    \label{fig:testcase}
\end{figure}
 
\section{EVCS Problem Formulation} \label{Problem Formulation}
The distribution grid and its EV aggregators are shown in \figref{fig:testcase}, which have a hierarchical trilayer structure including DSO, Aggs, and EVs. The goal of the DSO is to minimize the energy loss in the distribution grid lines, while Aggs aim at reducing the electricity cost for their customers. Accordingly, the objective function of the EV charging scheduling is twofold, loss reduction (grid service) and cost reduction (customer service), and it is written as follows:
\begin{equation} \label{ObjFun1}
\begin{split}
&V:=\min_{}\sum\limits_{(i,j) \in \zeta}^{}{r_{ij}\mathbf{I}_{ij}}+\sum\limits_{j\in \mathbb{N}_{Agg}}^{}{\mathbf{\Pi}^T.\mathbf{P}_j}\\
& \quad \text{s.t.}~\eqref{Grid2}-\eqref{Grid3},~\eqref{RESSet}-\eqref{Grid4},
\end{split}
\end{equation}
where $\mathbf{\Pi} \in \mathbb{R}^{N}$ is the wholesale energy price.
\subsection{Distributed EVCS}
As it is clear, solving the optimization problem in \eqref{ObjFun1} by only DSO is not computationally efficient, especially when the grid is bulky and DSO has to deal with a large population of EVs. Therefore, we use ADMM to solved the problem in a distributed manner in which DSO and Aggs communicate which each other iteratively. However, we first need to relax the optimization problem as \eqref{Grid24} is a non-convex constraint, otherwise ADMM is not applicable. According to \cite{Low1}, \eqref{Grid24} can be relaxed to a convex second-order cone as follows:  
\begin{equation}
{\mathbf{V}}_i{\mathbf{I}}_{ij}\geq\mathbf{P}_{ij}^2+\mathbf{Q}_{ij}^2.\label{RelGrid24}
\end{equation}
The sufficient conditions to make the relaxation tight, as shown in \cite{Low1}-\cite{Low2}, are: (i) the grid should be radial; (ii) bus voltages should be very close to the nominal value; and (iii) the active and reactive powers injected to the buses should not be too large. 

Replacing \eqref{Grid24} by \eqref{RelGrid24} in \eqref{ObjFun1}, the optimization problem \eqref{ObjFun1} is solved by ADMM as the following:
\begin{align} \label{ADMM1_x1}
\begin{split} 
&(\mathbf{P}_{j_c}^{k+1},\mathbf{Q}_{j_c}^{k+1}):=\argmin_{\mathbf{P}_{j_c},\mathbf{Q}_{j_c}} \Big(\mathbf{\Pi}^T.\mathbf{P}_{j_c}+\frac{\rho_p}{2}\norm{\mathbf{P}_{j_c} - \mathbf{P}_{j}^{k} + \mathbf{v}_{j}^{k}}_{2}^{2}\\ & +\frac{\rho_q}{2}\norm{\mathbf{Q}_{j_c} - \mathbf{Q}_{j}^{k} + \mathbf{u}_{j}^{k}}_{2}^{2} \Big) \\
& \qquad \qquad \text{s.t.}~\mathbb{U}_{j,i}~\forall i \in \mathbb{N}_{j,ev}~\&~j \in \mathbb{N}_{Agg},~\eqref{BES}-\eqref{Grid4}
\end{split}
\end{align}
\begin{align} \label{ADMM1_y1}
\begin{split}
&(\mathbf{P}_{j}^{k+1},\mathbf{Q}_{j}^{k+1}):=\argmin_{\mathbf{P},\mathbf{Q},\mathbf{V},\mathbf{I}} \Big(\sum\limits_{(i,j) \in \zeta}^{}{r_{ij}\mathbf{I}_{ij}}\\
&+\frac{\rho_p}{2}\sum\limits_{j\in \mathbb{N}_{Agg}}^{}{\norm{\mathbf{P}_{j_c}^{k+1} - \mathbf{P}_{j} + \mathbf{v}_{j}^{k}}_{2}^{2}}\\ 
& +\frac{\rho_q}{2}\sum\limits_{j\in \mathbb{N}_{Agg}}^{}{\norm{\mathbf{Q}_{j_c}^{k+1} - \mathbf{Q}_{j} + \mathbf{u}_{j}^{k}}_{2}^{2}} \Big) \\
&\quad \quad \text{s.t.}~\eqref{Grid21}-\eqref{Grid23},~\eqref{Grid3}~\& ~\eqref{RelGrid24}
\end{split}
\end{align}
\begin{subequations} \label{ADMM1_u1}
\begin{align}
&\mathbf{v}_{j}^{k+1} = \mathbf{v}_{j}^{k}+\mathbf{P}_{j_c}^{k+1}-\mathbf{P}_{j}^{k+1} \label{ADMM1_u11} \\ 
&\mathbf{u}_{j}^{k+1} = \mathbf{u}_{j}^{k}+\mathbf{Q}_{j_c}^{k+1}-\mathbf{Q}_{j}^{k+1}, \label{ADMM1_u12}
\end{align}
\end{subequations}
where \eqref{ADMM1_x1} is solved in parallel by each Agg, and \eqref{ADMM1_y1}-\eqref{ADMM1_u1} are solved by DSO.
\subsection{Hierarchical Distributed EVCS}
Considering the first step of ADMM \eqref{ADMM1_x1}, each Agg has to solve the optimal charging scheduling problem for all the EVs which it is supplying. If the number of EVs is considerable, the computational burden for Aggs will be substantial. By exploiting the mathematical formulation, we write \eqref{ADMM1_x1} in the form which is called \textit{sharing problem}, and it can be solved efficiently by ADMM in a distributed manner between each $Agg_j$ and its EVs (i.e. $\forall EV_{j,i}$, $i \in \mathbb{N}_{j,ev}$).

Using \eqref{Grid4} and \eqref{Grid5}, we can write \eqref{ADMM1_x1} as:
\begin{align} \label{ADMM1_x2}
\begin{split} 
&\min_{x,y}(\sum\limits_{i \in \mathbb{N}_{j,ev}}^{}{\mathbf{\Pi}^T.\mathbf{p}_{j,i}^{ev}})+\mathbf{\Pi}^T.(\mathbf{p}_{j}^{uc}+\mathbf{p}_{j}^{bes})\\
&+\frac{\rho_p}{2}\norm{\sum\limits_{i \in \mathbb{N}_{j,ev}}^{}{\mathbf{p}_{j,i}^{ev}}+ \mathbf{p}_{j}^{uc}+\mathbf{p}_{j}^{bes}- \mathbf{P}_{j_c}^{k} + \mathbf{v}_{j}^{k}}_{2}^{2}\\ 
& +\frac{\rho_q}{2}\norm{\mathbf{q}_{j}^{uc}+\mathbf{q}_{j}^{bes} - \mathbf{Q}_{j_c}^{k} + \mathbf{u}_{j}^{k}}_{2}^{2} \\
& \qquad \text{s.t.}~\mathbb{U}_{j,i}~\forall i \in \mathbb{N}_{j,ev}~\&~j \in \mathbb{N}_{Agg},~\eqref{BES}-\eqref{Grid4}.
\end{split}
\end{align}
Considering $\mathbf{p}_{j,i}^{ev}$, the first part on RHS of \eqref{ADMM1_x2} is separable between the EVs, and the rest is a function of $\sum\limits_{i \in \mathbb{N}_{j,ev}}^{}{\mathbf{p}_{j,i}^{ev}}$ which we show, hereafter, by $\mathbf{p}_{j}^{ev}$. Therefore, \eqref{ADMM1_x2} is a \textit{sharing problem} \cite[Chapter~7.3]{Boyd1}, and it can be solved in a distributed manner as follows:
\begin{equation} \label{ADMM2_x1}
\begin{split} 
&\mathbf{p}_{j,i}^{{ev}^{l+1}}:=\argmin_{\mathbf{p}_{j,i}^{ev}} \Big(\mathbf{\Pi}^T.\mathbf{p}_{j,i}^{ev}+\\
&\frac{\rho_{j}}{2}\norm{\mathbf{p}_{j,i}^{ev} - \mathbf{p}_{j,i}^{{ev}^{l}} + \overline{\mathbf{p}}_{j}^{{ev}^l}-\overline{\mathbf{p}}_{j_{c}}^{{ev}^l}+\pmb{\lambda}_j^{k}}_{2}^{2}\Big)\\
&\qquad \qquad \text{s.t.}~\mathbb{U}_{j,i}~\forall i \in \mathbb{N}_{j,ev}
\end{split}
\end{equation}
\begin{align} \label{ADMM2_y1}
\begin{split}
&\overline{\mathbf{p}}_{j_{c}}^{{ev}^{l+1}}:=\argmin_{\mathbf{p}_{j}^{bes},\mathbf{q}_{j}^{bes},\overline{\mathbf{p}}_{j_{c}}^{ev}}\mathbf{\Pi}^T.(\mathbf{p}_{j}^{uc}+\mathbf{p}_{j}^{{bes}})\\
&+\frac{\rho_p}{2}\norm{\mathcal{N}_{j,ev}.\overline{\mathbf{p}}_{j_{c}}^{ev}+ \mathbf{p}_{j}^{uc}+\mathbf{p}_{j}^{{bes}}- \mathbf{P}_{j_c}^{k} + \mathbf{v}_{j}^{k}}_{2}^{2}\\ 
& +\frac{\rho_q}{2}\norm{\mathbf{q}_{j}^{uc}+\mathbf{q}_{j}^{{bes}} - \mathbf{Q}_{j_c}^{k} + \mathbf{u}_{j}^{k}}_{2}^{2} \\
& +(\frac{\mathcal{N}_{j,ev}.\rho_{j}}{2})\norm{\overline{\mathbf{p}}_{j_{c}}^{ev}-\overline{\mathbf{p}}_{j}^{{ev}^{l+1}}-\pmb{\lambda}_j^{k}}_{2}^{2}\\
& \qquad\qquad\qquad\qquad \text{s.t.}~\eqref{BES}-\eqref{BESConstraint}
\end{split}
\end{align}
\begin{equation} \label{ADMM2_u1}
\pmb{\lambda}_{j}^{l+1} = \pmb{\lambda}_{j}^{l}+\overline{\mathbf{p}}_{j}^{{ev}^{l+1}}-\overline{\mathbf{p}}_{j_{c}}^{{ev}^{l+1}}. 
\end{equation}
Note that $Agg_j$'s problem size \eqref{ADMM2_y1} is independent of the number of EVs. To decrease communication overheads, we define $\pmb{\Lambda}_j^{l+1}=\pmb{\lambda}_{j}^{l+1}+\overline{\mathbf{p}}_{j}^{{ev}^{l+1}}-\overline{\mathbf{p}}_{j_{c}}^{{ev}^{l+1}}$. Thus, at each \textit{sharing problem} iteration after the third update \eqref{ADMM2_u1}, $Agg_j$ broadcasts $\pmb{\Lambda}_j^{l+1}$ to all $EV_{j,i}$, $\forall i \in \mathbb{N}_{j,ev}$. For more details, we refer the readers to \cite{Behnam1}. 

We call \textit{ADMM}\textsubscript{1} the iterative procedure between DSO and Aggs \eqref{ADMM1_x1}-\eqref{ADMM1_u1}, and \textit{ADMM}\textsubscript{2} the \textit{sharing problem} between each Agg and its EVs \eqref{ADMM2_x1}-\eqref{ADMM2_u1}. The whole procedure of our hierarchical distributed EVCS is shown in Algorithm \ref{algorithm1}. $Err_1$ and $Err_2$ are the pair of primal and dual residuals for \textit{ADMM}\textsubscript{1} and \textit{ADMM}\textsubscript{2}, respectively. For more details about residual calculation and stopping criteria, we refer the reader to \cite[Chapter~3.3]{Boyd1} 
\begin{algorithm}
\SetAlgoLined
\DontPrintSemicolon
%\KwData{Calculate $\overline{\Omega}(t)$}
%\KwResult{Write here the result }
% \textbf{Initialization} Calculate $\overline{\Omega}(t)$ \& EVs' arrival/departure time and charging energy \;
 \While{$Err_{1}<Th_{1}$}{
      \For{$j \in \mathbb{N}_{Agg}$}{
       \While{$Err_{2}<Th_{2}$}{
        \For{$i \in \mathbb{N}_{j,ev}$}{
         Calculate $\mathbf{p}_{j,i}^{{ev}}$ by \eqref{ADMM2_x1} \& send to $Agg_j$.
         }
        Update $\overline{\mathbf{p}}_{j}^{ev}=\frac{1}{\mathcal{N}_{j,ev}}\sum\limits_{i \in \mathbb{N}_{j,ev}}^{}\mathbf{p}_{j,i}^{{ev}}$.\;
        Calculate $\overline{\mathbf{p}}_{j_{c}}^{ev}$, $\mathbf{p}_{j}^{{bes}^{l+1}}$ and $\mathbf{q}_{j}^{{bes}^{l+1}}$ by \eqref{ADMM2_y1}.\;
        Update $\pmb{\lambda}_{j}$ by \eqref{ADMM2_u1}.\;
        Update $\pmb{\Lambda}_{j}$ and broadcast to $\forall i \in \mathbb{N}_{j,ev}$.\;
        Update $Err_{2}$.\;
       }
       Send $(\mathbf{P}_{j},\mathbf{Q}_{j})$ to DSO.
      }
   Calculate $(\mathbf{P}_{j_c},\mathbf{Q}_{j_c})$, $\forall j \in \mathbb{N}_{Agg}$, by \eqref{ADMM1_y1}.\;
   Update $(\mathbf{v}_{j},\mathbf{u}_{j})$, $\forall j \in \mathbb{N}_{Agg}$, by \eqref{ADMM1_u1}.\;
   Broadcast $(\mathbf{P}_{j_c},\mathbf{Q}_{j_c})$ and $(\mathbf{v}_{j},\mathbf{u}_{j})$ to $Agg_j$, $\forall j \in \mathbb{N}_{Agg}$. 
   Update $Err_{1}$.\;
 }
\caption{Hierarchical Distributed EVCS.}
\label{algorithm1}
\end{algorithm}

\section{Numerical Simulation} \label{Simulation Results}
In this section, the performance of the proposed hierarchical distributed EVCS is evaluated for the modified IEEE-13 bus system. We consider the signle phase balanced system with six Aggs which are located at bus\# $634$, $646$, $675$, $680$, $652$ and $611$. We compare the performance of our proposed method with uncontrolled EV charging, in which EVs start charging as soon as they are plugged in with the maximum power rating (i.e. $\overline{\mathbf{p}}^{ev}_{j,i}$, $\forall i \in \mathbb{N}_{j,ev}$, $\forall j \in \mathbb{N}_{Agg}$). Also, to show the effect of BES on loss and charging cost reduction, we compare the results with the case where EVCS is executed without any BES installed in the Aggs' nodes. The maximum power rating for all EV chargers is $4$ kW. The initial and designated EVs' battery energies are uniformly distributed over $[8,10]$ kWh and $[22,25]$ kWh, respectively.  Also, EVs' arrival and departure times are normally distributed in $[\texttt{16:30},\texttt{20:30}]$ and $[\texttt{6:00},\texttt{9:30}]$, respectively. The netload dataset of EVBs is collected from the Australian electricity company-Ausgrid \cite{Ausgrid}, and the wholesale price is available from the California Independent System Operator-CAISO \cite{CAISO}. More details about the simulation parameters are shown in Table \ref{SimParameters}. All the simulations are executed by MATLAB on a PC with Intel$\text{\textregistered}$ Core$\text{\texttrademark}$ i7-4770 3.40 GHz CPU, 4 cores and 8 GB RAM, and the convex optimization problems are solved by CVX \cite{CVX}.
\begin{table}[h]
\caption{EV, BES and hierarchical EVCS simulation parameters}\label{SimParameters}
\centering\setlength{\arrayrulewidth}{0.9pt}
%\begin{tabular}{@{} ll|ll @{}}
\begin{tabular}{ cc | cc }
%\toprule
\hline
%\specialrule{.14em}{.08em}{.05em} 
Parameter & Value & Parameter & Value \\
%& (s) & (kbps) & \\
%\midrule
%\specialrule{.15em}{.1em}{.1em} 
\hline
$\overline{\mathbf{p}}^{ev}$, $\underline{\mathbf{p}}^{ev}$ & $4$, $-4$ [kW] &$\mathbf{v}_{ref}$ & 1.0 [p.u.]\\
$\overline{\mathbf{v}}$, $\underline{\mathbf{v}}$& 0.97, 1.03   [p.u.] &$\rho_p$, $\rho_q$ & 0.1, 0.1 \\
$\overline{\mathbf{C}}_{bes}$ & $55$ [kWh] & $\rho_j$  & 1 \\
$\underline{\mathbf{C}}_{bes}$& $5$ [kWh] &$N$ & $48$ \\
$\overline{\mathbf{s}}_{bes}$ & $50$ [kVA]& $T_h $ & $0.5$ \\
%\bottomrule
\hline
\end{tabular}
\end{table}

\begin{figure}[h!]\hspace*{-.2cm}
	\centering
	\subfloat %[\acs{SD}.]
	{\includegraphics[clip,scale=0.5, trim=1.2cm 9.5cm 14.4cm 11.6cm]{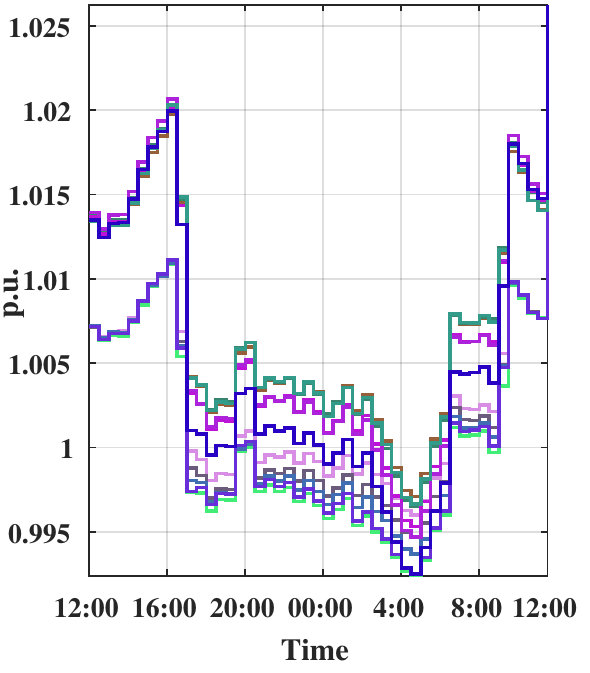}\label{VolwBES}}%\hfill
	\subfloat %[\acs{DE}.]
	{\includegraphics[clip,scale=0.5, trim=1.1cm 9.5cm 14.4cm 11.6cm]{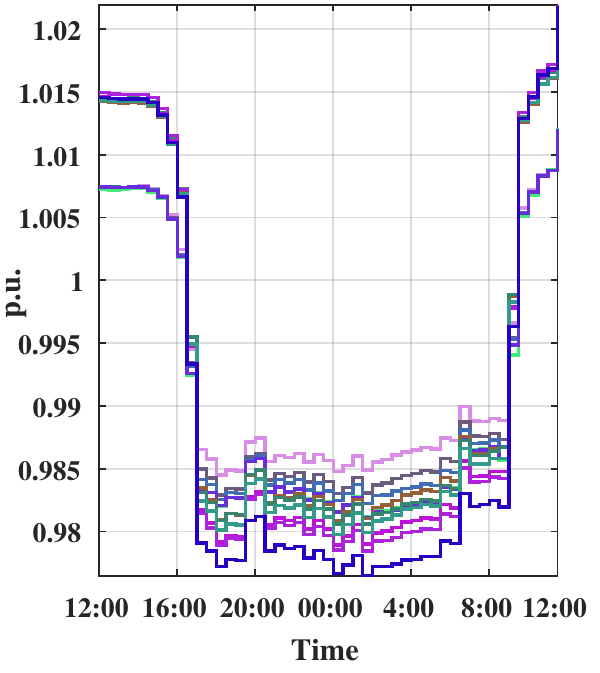}\label{VolwoBES}}
		\subfloat %[\acs{DE}.]
	{\includegraphics[clip,scale=0.5, trim=1.2cm 9.5cm 14.0cm 11.6cm]{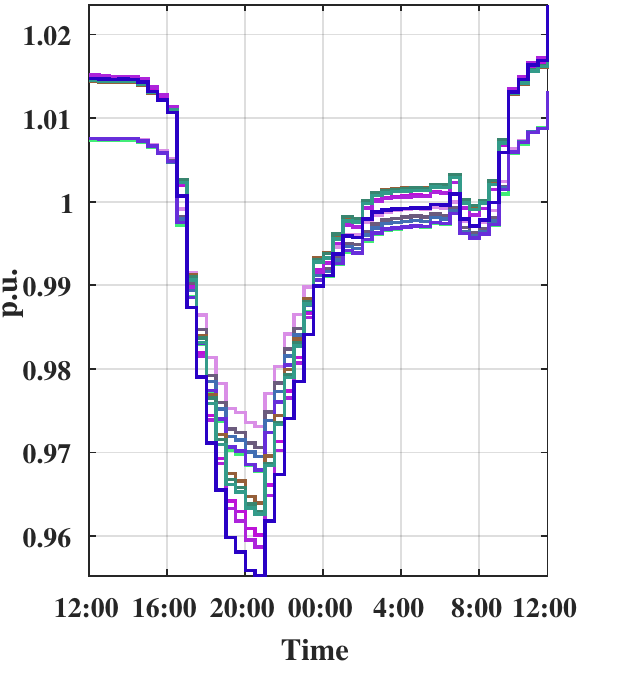}\label{VolwoBES}}
	\caption{Voltage profiles of the grid buses: (left) CC\textsubscript{1}, (middle) CC\textsubscript{2}, and (right) uCC.}\label{Vol}
\end{figure}

As maintaining the bus voltages within the acceptable range is a constraint of EVCS, \figref{Vol} shows how the control EV charging demand with (CC\textsubscript{1}) and without stationary BES (CC\textsubscript{2}) improves voltage profile in the grid buses while uncontrolled EV charging (uCC) results in significant voltage drop in the grid. The results also show that EVCS with BES results in better voltage profile, i.e. bus voltages are closer to the nominal value ($1.0$ [p.u.]). 
%\begin{figure}[h!]\hspace*{-.2cm}
%	\centering
%	\subfloat %[\acs{SD}.]
%	{\includegraphics[clip,scale=0.5, trim=0.6cm 9.4cm 12.0cm 11.5cm]{VoltageProfile_03.pdf}\label{VolwBES}}%\hfill
%	\subfloat %[\acs{DE}.]
%	{\includegraphics[clip,scale=0.5, trim=0.6cm 9.4cm 12.0cm 11.5cm]{VoltageProfile2_01.pdf}\label{VolwoBES}}\\
%	\caption{Voltage profiles of the grid buses, (left) EVCS with BES, (right) EVCS without BES: uncontrolled charging (dashed line), proposed EVCS (solid line).}\label{Vol}
%\end{figure}
\begin{figure}[h!]\hspace*{-.2cm}
	\centering
	\subfloat %[\acs{SD}.]
	{\includegraphics[clip,scale=0.5, trim=1.1cm 9.4cm 11.7cm 11.5cm]{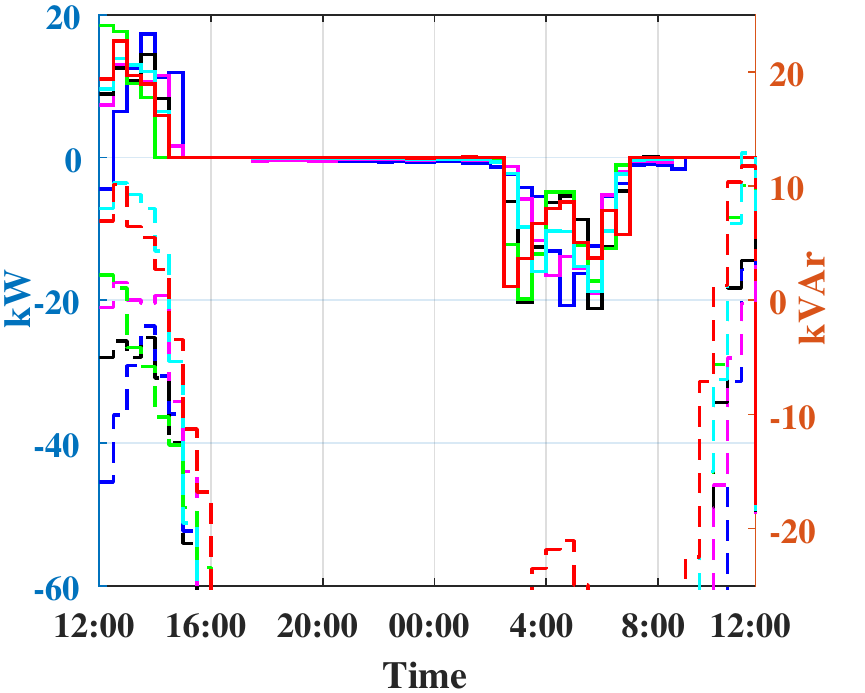}\label{BESPQ}}%\hfill
	\subfloat %[\acs{DE}.]
	{\includegraphics[clip,scale=0.5, trim=1.2cm 9.4cm 11.4cm 11.5cm]{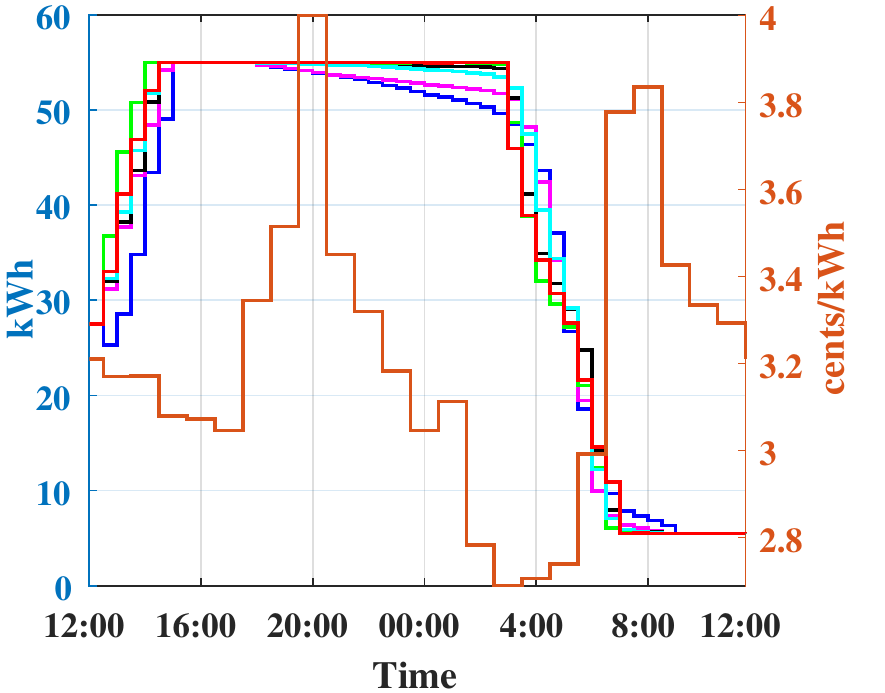}\label{BESTOU}}\\
	\caption{(left) BES active and reactive powers, (right) BES energy profile and wholesale electricity price.}\label{BESPQETOU}
\end{figure}

\begin{figure}[h!]\hspace*{-.2cm}
	\centering
	\subfloat %[\acs{SD}.]
	{\includegraphics[clip,scale=0.5, trim=1.6cm 9.5cm 14.1cm 11.6cm]{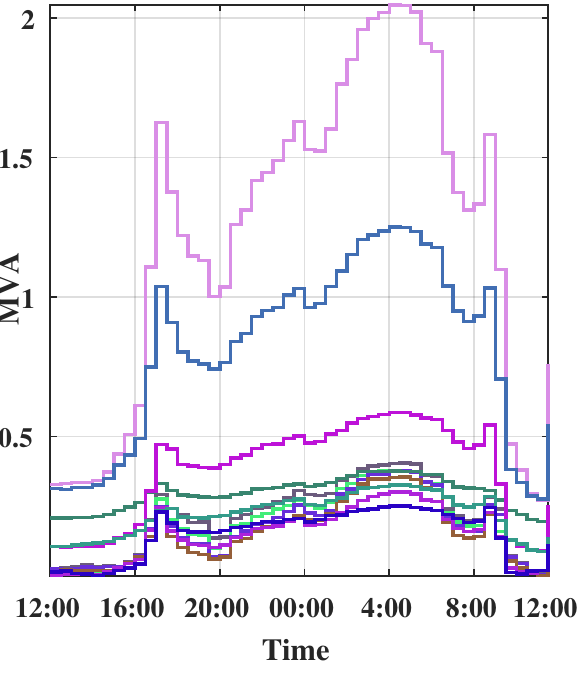}\label{LPCC1}}%\hfill
	\subfloat %[\acs{DE}.]
	{\includegraphics[clip,scale=0.5, trim=1.6cm 9.5cm 14.1cm 11.6cm]{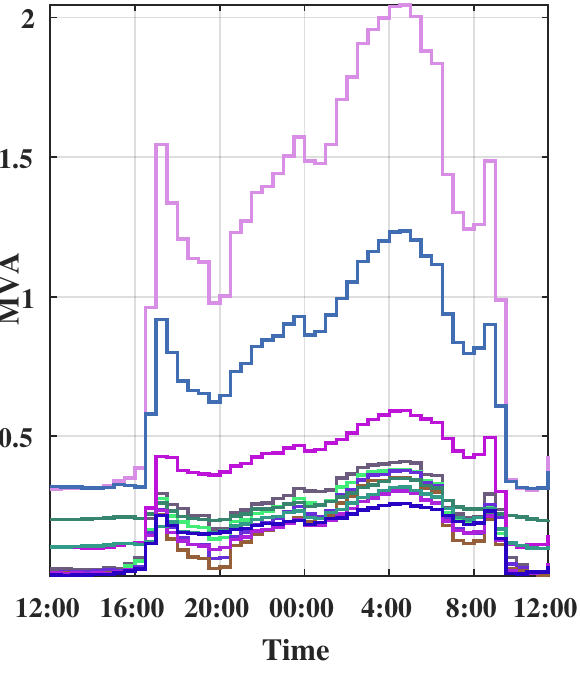}\label{LPCC2}}
		\subfloat %[\acs{DE}.]
	{\includegraphics[clip,scale=0.5, trim=1.3cm 9.5cm 14.0cm 11.6cm]{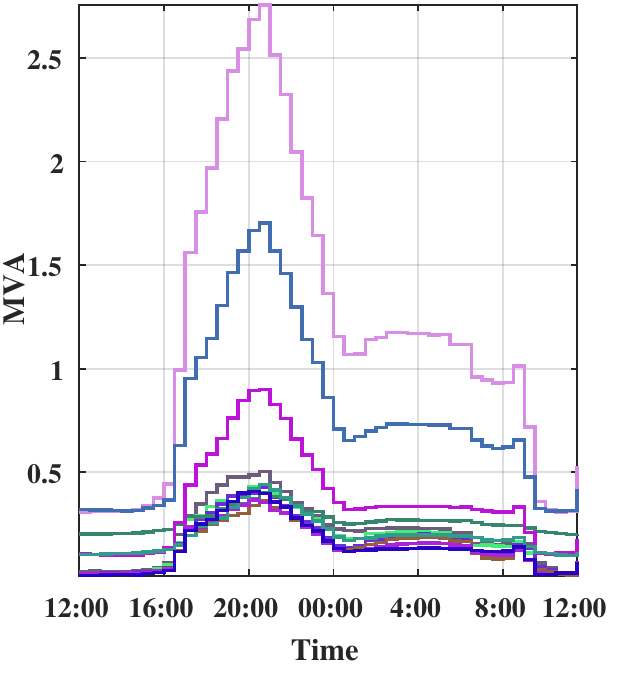}\label{LPuCC}}
	\caption{Load profile of the grid lines: (left) CC\textsubscript{1}, (middle) CC\textsubscript{2}, and (right) uCC.}\label{LProfile}
\end{figure}

\figref{BESPQETOU} shows the active and reactive powers and the energy profile of the BES as well as the wholesale electricity price. As it is shown, BES is charged while the energy price is low, and it is discharged during the second peak of the electricity price. It should be noticed that BES is not discharged during the first price peak as load demand is not considerable (\figref{LProfile}). While close to the second price peak when the load demand is considerable, BES is discharged to reduce both energy loss and electricity cost. Considering \figref{LProfile}, the peak load over the incoming transformer feeder decreases considerably in both CC\textsubscript{1} and CC\textsubscript{2} ($\sim 2$ MVA), while its is $\sim 2.7$ MVA for uCC.  

\begin{figure}[h]\hspace*{0cm}
    \centering \includegraphics[clip,scale=0.53, trim=1.15cm 9.8cm 3.cm 12.7cm] %width=0.5\textwidth
    {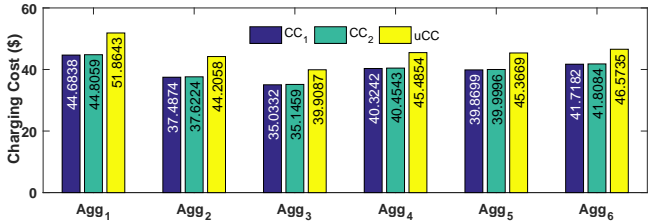}%{test255075-2.pdf}
        \caption{{The comparison of EV charging costs obtained by CC\textsubscript{1}, CC\textsubscript{2}, and uCC.}}
%    \caption{PnP using the proposed MPC-ADMM; switching between charging cost and battery degradation reduction modes}
    \label{fig:AggCost}
\end{figure}

In \figref{fig:AggCost}, the EV charging costs using the three simulated methods are compared. As it is expected, CC\textsubscript{1} achieves the least cost, and uCC results in the highest cost.
%\figref{fig:AggCost} shows the comparison between the aggregated EV charging cost using the three simulated methods. As it is expected, CC\textsubscript{1} propose the least cost, and uCC results in the highest cost.

In addition, DSO achieves its energy loss reduction in the grid as it is shown in Table \ref{LossRed}. The loss reduction achieved by CC\textsubscript{1} and CC\textsubscript{2} is considerable in the lines\# $1$, $4$ and $8$, which connect bus\# $650$ to $632$, $632$ to $645$, and $671$ to $684$, respectively. The reasons for the significant loss are high impedance of the line conductors, their length, and high peak loads.
\begin{table}[h!]
	\caption{Energy loss on the grid lines}\label{LossRed}
	\centering
	   \begin{threeparttable}[b]
	\setlength\tabcolsep{3pt}
	\scalebox{0.91}{
	\begin{tabular}{@{} ccccccccccccc @{}}
		\hlinewd{1pt}
		\addlinespace[0.05cm]
		 {}& Loss & 1 & 2 & 3 & 4 & 5 & 6 & 7 & 8 & 9 & 10 & 11 \\ \hlinewd{1pt}
		\multirow{3}{*}{{kWh}} &{CC\textsubscript{1}} & $325$ & $6.11$ & $8.65$ & $117$  & $6.46$ & $3.50$ & $7.36$ & $16.60$ & $2.60$ & $4.60$ & $8.08$\\
		&{CC\textsubscript{2}} & $331$ & $6.40$ & $9.08$ & $115$  & $6.67$ & $3.59$ & $6.68$ & $16.75$ & $2.71$ & $4.34$ & $8.38$\\
		&{uCC}  & $375$ & $6.86$ & $9.34$ & $133$  & $6.76$ & $3.65$ & $6.92$ & $20.2$ & $2.99$ & $4.95$ & $10.6$\\ \hlinewd{0.5pt}
		\addlinespace[0.05cm]
		\multirow{3}{*}{{kVarh}} & {CC\textsubscript{1}} & $102$ & $9.80$ & $8.70$ & $35.4$  & $1.20$ & $3.60$ & $4.10$ & $17.1$ & $7.70$ & $4.70$ & $3.10$\\
		& {CC\textsubscript{2}} & $107$ & $10.3$ & $9.30$ & $34.9$  & $1.20$ & $3.70$ & $3.70$ & $17.2$ & $8.10$ & $4.40$ & $3.20$\\
		&{uCC} & $123$ & $11.0$ & $9.60$ & $40.15$  & $1.20$ & $3.40$ & $3.90$ & $20.7$ & $8.90$ & $5.21$ & $4.03$\\
		\addlinespace[0.05cm]
		\hlinewd{1pt}
	\end{tabular}}
%	   \begin{tablenotes}
%     \item[1] controlled charging by the hierarchical EVCS with BES.
%     \item[2] controlled charging by the hierarchical EVCS without BES.
%     \item[3] uncontrolled charging.
%   \end{tablenotes}
  \end{threeparttable}
\end{table}
%As \textit{EVB}s and \textit{Agg}s try to minimize the electricity cost, EVCS achieve this goal by optimal charging/discharging of the \textit{EV}s and BESs, as it is shown in \figref{ChCost}. According to the twofold objective function of EVCS, there is a trade-off between energy loss reduction and electricity cost reduction. By defining weighting factors for each objective, the extend to which energy loss or cost is reduced varies.  
\section{Conclusion} \label{Conclusion}
In this paper, a fully hierarchical and distributed method was proposed for EVCS. While the power flow model and constraints are considered for optimal scheduling of EV load demand, DSO in collaboration with Aggs and EVs try to minimize the energy loss and electricity cost. As solving EVCS optimization problem in a centralized manner is not computationally efficient and privacy preserving, ADMM was used to solve the problem through an iterative procedure between DSO, Aggs, and EVs. Numerical simulation of the proposed EVCS, which was applied to IEEE-13 bus system including six Aggs, verified its effectiveness. 
\section*{Acknowledgment}
The first author would like to acknowledge American Public Power Association (APPA) for supporting his research through the DEED grant.
%\begin{figure}[h!]\hspace*{-.2cm}
%	\centering
%	\subfloat %[\acs{SD}.]
%	{\includegraphics[clip,scale=0.5, trim=1.2cm 9.5cm 14.4cm 11.6cm]{NodVoltage_CC1_01.pdf}\label{VolwBES}}%\hfill
%	\subfloat %[\acs{DE}.]
%	{\includegraphics[clip,scale=0.5, trim=1.1cm 9.5cm 14.4cm 11.6cm]{NodVoltage_CC2_01.pdf}\label{VolwoBES}}
%		\subfloat %[\acs{DE}.]
%	{\includegraphics[clip,scale=0.5, trim=1.2cm 9.5cm 14.0cm 11.6cm]{NodVoltage_nCC_01.pdf}\label{VolwoBES}}
%	\caption{Voltage profiles of the grid buses, (left) EVCS with BES, (right) EVCS without BES: uncontrolled charging (dashed line), proposed EVCS (solid line).}\label{Vol}
%\end{figure}

\end{document}